\documentclass{article}
\usepackage[utf8]{inputenc}
\usepackage{amsmath}
\usepackage{amssymb}
\usepackage{amsthm}
\usepackage{verbatim}
\usepackage{graphicx}
\usepackage{subcaption}
\usepackage[shortlabels]{enumitem}
\usepackage{xcolor}
\usepackage{float}
\usepackage{lineno}
\usepackage{cases}
\usepackage{bm}
\usepackage[normalem]{ulem}

\newtheorem{theorem}{Theorem}[section]
\newtheorem{corollary}{Corollary}[theorem]
\newtheorem{lemma}[theorem]{Lemma}

\newtheorem{remark}[theorem]{Remark}

\newcommand{\aam}{{a_{j}}}

\newcommand{\gamt}{{\tilde{\gamma}_{j}}}
\newcommand{\notj}{{\tilde{j}}}

\newcommand{\FhD}{{F_{h,\Diamond}}}
\newcommand{\AAA}{{\mathsf{A}}}

\newcommand{\AAd}{{\mathsf{A}^\dagger}}
\newcommand{\NAT}{{\mathcal{N}(\AAA)^\bot}}

\newcommand{\JJ}{\mathcal{J}}
\newcommand{\sJ}{{\sum_{\mathcal{J}}}}
\newcommand{\sJee}{{\sum_{\mathcal{J}}x_j^*\ee_j}}

\newcommand{\sJgtee}{{\sum_{\mathcal{J}}\tilde{\gamma}_j\ee_j}}

\newcommand{\WW}{{\mathsf{W}}}

\newcommand{\PP}{\mathsf{P}}

\newcommand{\II}{{\mathsf{I}}}
\newcommand{\xx}{{\mathbf{x}}}
\newcommand{\yy}{{\mathbf{y}}}
\newcommand{\zz}{{\mathbf{z}}}
\newcommand{\vv}{{\mathbf{v}}}

\newcommand{\bb}{{\mathbf{b}}}
\newcommand{\cc}{{\mathbf{c}}}

\newcommand{\qq}{{\mathbf{q}}}
\newcommand{\ee}{{\mathbf{e}}}
\newcommand{\rr}{{\mathbf{r}}}

\newcommand{\Rn}{{\mathbb{R}^n}}

\newcommand{\Rnn}{{\mathbb{R}^{n\times n}}}
\newcommand{\Rmn}{{\mathbb{R}^{m\times n}}}

\newcommand{\nullspace}[1]{{\mathcal{N}(#1)}}

\DeclareMathOperator*{\argmin}{arg\,min}
\DeclareMathOperator*{\argmax}{arg\,max}

\newcommand{\rem}[1]{\textcolor{black}{#1}}

\title{Identifying the source term in the potential equation with weighted sparsity regularization}

\author{Ole L{\o}seth Elvetun\thanks{Faculty of Science and Technology, Norwegian University of Life Sciences, P.O. Box 5003, NO-1432 {\AA}s, Norway. Email: ole.elvetun@nmbu.no.} and Bj{\o}rn Fredrik Nielsen\thanks{Faculty of Science and Technology, Norwegian University of Life Sciences, P.O. Box 5003, NO-1432 {\AA}s, Norway. Email: bjorn.f.nielsen@nmbu.no.}}

\begin{document}

\maketitle

\begin{abstract}
  We explore the possibility for using boundary measurements to recover a sparse source term $f(x)$ in the potential equation. Employing weighted sparsity regularization and standard results for subgradients, we derive simple-to-check criteria which assure that a number of sinks ($f(x)<0$) and sources ($f(x)>0$) can be identified. 
  Furthermore, we present two cases for which these criteria always are fulfilled: a) well-separated sources and sinks, and b) many sources or sinks located at the boundary plus one interior source/sink. 
  Our approach is such that the linearity of the associated forward operator is preserved in the discrete formulation. The theory is therefore conveniently developed in terms of Euclidean spaces, and it can be applied to a wide range of problems. In particular, it can be applied to both  isotropic and anisotropic cases.  We present a series of numerical experiments. 
  This work is motivated by the observation that standard methods typically suggest that internal sinks and sources are located close to the boundary.
\end{abstract}

\noindent {\bf 2020 Mathematics Subject Classification:} 35R30, 47A52, 65F22. \\

\noindent {\bf Keywords:}
Inverse source problems, potential equation, anisotropy, null space, sparsity regularization, PDE-constrained optimization.

\section{Introduction}
The boundary value problem 
\begin{align}
    -\nabla \cdot \sigma\nabla u &= f, \quad x \in \Omega, \label{eq:in1} \\
    \mathbf{n}\cdot \sigma\nabla u &= 0, \quad x \in \partial\Omega,\label{eq:in2} \\ 
    \int_{\partial\Omega} u &= 0, \label{eq:in3}
\end{align}
appears in several applications: voltage models, heat conduction, descriptions of gravitational fields, flow in porous media, etc.  
In such models, $u$ represents the state, $\sigma=\sigma(x)$ is a conductivity tensor or function, $f$ is the source term,
equation \eqref{eq:in2} assures insulation at $\partial \Omega$ and \eqref{eq:in3} is included to enforce uniqueness of the solution $u$ of \eqref{eq:in1}-\eqref{eq:in2}. \rem{Throughout this paper we assume that $\partial \Omega$ is piecewise smooth.} 

Typically, the state function $u$, or some noisy version of it, is observed on the boundary $\partial\Omega$. Using such a recording $d$, the goal of the inverse source problem is to recover $f$ by solving  
\begin{equation} \label{eq:in3.5}
\min_{f,u} \| u -d \|_{L^2(\partial \Omega)} \quad \textnormal{subject to \eqref{eq:in1}-\eqref{eq:in3}}.    
\end{equation}
This problem is ill-posed: The associated forward operator has a nontrivial null space and its Moore-Penrose inverse is unbounded. 

The regions where, e.g., heat or charge is added or extracted may be small. 
In such cases, the sources ($f(x)>0$) and sinks ($f(x)<0$) will have small local support. 
Hence, we will in this paper employ weighted sparsity regularization to a discrete version of \eqref{eq:in3.5}. The weighting is needed because standard regularization methods fail to correctly identify sinks and sources which are far from the boundary (where the data is recorded), see, e.g., \rem{\cite{burger2013inverse,casas2012approximation,Elv21c,pascual2002standardized}.}

Problems similar to \eqref{eq:in3.5} have been studied by many mathematicians, especially the case with constant, or piecewise constant, conductivity $\sigma$. Often one assumes that $f$ is a sum of point sources/sinks, or dipoles, and one seeks to compute the positions of these localized functions and/or their strengths or moments, see, e.g., \cite{babda09,Chung_2009,Bad00}. This approach typically yields nonlinear systems of algebraic equations, but reliable and fast methods have been developed. It is also possible to determine the number of point sources with this type of methodology \cite{Bad00}. Furthermore, see \cite{Badia_2005} for results for the anisotropic case. 

The task of approximately determining the support of the right-hand-side $f$ in \eqref{eq:in1}-\eqref{eq:in3} from boundary data has also been studied in detail, see, e.g., \cite{Han11,Het96,BIsa05,liu2017inverse}. This leads to involved \rem{analysis} since the problem is ill posed and due to the large null space of the associated forward operator. Roughly speaking, one can not expect to determine accurate information about the position, the size and the magnitude of the sources and the sinks without imposing further apriori restrictions \cite{Han11}. 

If one employs the ansatz that $f$ is composed of a number of point sources and sinks, then the parameter-to-observation map, which, e.g., maps the positions of the point sources to the potential at the surface $\partial \Omega$, typically becomes nonlinear. This is in contrast to the linearity of the forward operator $f \mapsto u|_{\partial \Omega}$ associated with \eqref{eq:in1}-\eqref{eq:in3.5}. 

\rem{
As mentioned above, if the objective is to recover sources and sinks with small support, one could also attempt to formulate the problem using sparsity regularization. Applying \textit{true} sparsity (in a general finite-dimensional setting) would imply solving a problem in the form
\begin{equation*}
    \min_{\xx\in\Rn} \|\xx\|_0 \quad \textnormal{ subject to } \quad \AAA\xx = \bb, 
\end{equation*}
where the $\|\cdot\|_0$-"norm" simply counts the number of non-zero components in the vector $\xx$. This is, however, an NP-hard problem, and therefore one typically rather considers the "relaxed" convex problem
\begin{equation}
    \min_{\xx\in\Rn} \|\xx\|_1 \quad \textnormal{ subject to } \quad \AAA\xx = \bb.\label{eq:BS1}
\end{equation}
}

\rem{
Since the forward mapping associated with \eqref{eq:in1}-\eqref{eq:in3.5}, $f \mapsto u_{\partial\Omega}$, is a mapping from a $d$-dimensional domain $\Omega$ to the $d-1$-dimensional surface $\partial\Omega$, the inverse problem will be underdetermined. This is inherited by any standard discretization of the problem, and consequently, the associated forward matrix $\AAA$ will have dimensions $\AAA \in \Rmn$, where $m < n$. 
}

\rem{
Problems in the form \eqref{eq:BS1}, with $\AAA \in \Rmn, m < n$, have been extensively studied in the compressed sensing (CS) community. If the sparsest vector $\xx^*$ which satisfies $\AAA\xx = \bb$ has $s$ non-zero components, we have several exact recovery criteria in the CS-literature which guarantee that $\xx^*$ is a (unique) minimizer of \eqref{eq:BS1}. The most classical conditions require either \textit{low incoherence} of the matrix \cite{Donoho03}, the fulfilment of the \textit{restricted isometry property} (RIP) \cite{candes05} or a certain bound on the \textit{exact recovery coefficient} (ERC) \cite{Fuchs04b}.
}

\rem{
Unfortunately, as we demonstrate in the next section, approximately solving \eqref{eq:BS1}, with the forward matrix $\AAA$ associated with our inverse problems, fails to recover even a 2-sparse solution. 
This is linked to the intrinsic depth bias of inverse source problems for elliptic PDEs, see, e.g., \cite{burger2013inverse,casas2012approximation,Elv22,lucka2012hierarchical}. That is, using standard methods, the boundary data will be "explained" with sources and sinks located (close to) where the data is observed. These observations provide both a computational verification that the involved matrix fails to meet the classical recovery criteria, mentioned in the previous paragraph, and that the use of $\ell_1$-regularization itself does not necessarily seem to give exact recovery for the present inverse source problem.
}

\rem{
In this paper we propose to overcome the previously mentioned depth bias issue by employing a tailored diagonal weight matrix and combine it with sparsity regularization.} This weighting procedure was also studied in \cite{Elv22} for the case when apriori box-constraints $0 \leq f(x) \leq \rem{f_{\max}}$, $x \in \Omega$, are available and when the state equation is either the screened Poisson equation or the Helmholtz equation. Here, $\rem{f_{\max}}$ is a given upper bound. However, in applications one will typically have both sources and sinks, 
and the right-hand-side $f$ in \eqref{eq:in1}-\eqref{eq:in3} must satisfy the complimentary condition   
\begin{equation} \label{eq:in3.6}
    \int_\Omega f = 0. 
\end{equation}
Also, there is no zero order term present in \eqref{eq:in1}. The results published in \cite{Elv22} can therefore not be applied, and a separate investigation is needed. 

Removing the box-constraints, and not requiring that $f$ is non-negative, makes the source-sink identification "harder". Nevertheless, we will in this paper prove that a number of well-separated sources and sinks can be recovered. It also turns out that an arbitrary number of sources or sinks located at the boundary $\partial \Omega$ plus one interior source/sink can be identified. Our analysis addresses both the basis pursuit version (section \ref{sec:basis_pursuit}) and the regularized version (section \ref{sec:regularized_problems}) of the problem. We also prove that a source will not be misinterpreted as a sink, or vice versa, and we present results regarding the uniqueness of the solution of the inverse problem. The analysis is complemented with numerical experiments in section \ref{sec:numerical_experiments}. 

One may regard this paper to be follow-up work to \cite{Elv21c}: In \cite{Elv21c} the single-source-case is analyzed when the state equation has the form $\Delta u + \epsilon u$, $\epsilon \neq 0$, and in all the numerical experiments in that paper the true source function is non-negative, $f(x) \geq 0$, $x \in \Omega$. 

We only consider finite dimensional problems in this paper, and some remarks concerning the discretization of \eqref{eq:in3.5} is presented in section \ref{sec:discretization}. In fact, our analysis is presented in terms of Euclidean spaces and can therefore be applied in other contexts, not just to the present inverse problem. 

 \rem{The idea of employing weighting procedures to the kind of inverse problem considered in this paper is not new. For example, this approach is addressed on pages 159--161 in \cite{burger2013inverse}, and the authors of that paper conclude that {\em "..., one observes that all standard approaches to incorporate
prior knowledge suffer from severe shortcomings in the application to imaging
from (underdetermined) surface data. It remains an important future challenge to develop
improved approaches that can provably reconstruct structures corresponding
to the available prior knowledge.".} Our aim is to contribute to solving this problem.}

\section{\rem{Motivation}}
\label{sec:motivation}
\begin{figure}[H]
    \centering
    \begin{subfigure}[b]{0.49\linewidth}        
        \centering
        \includegraphics[width=\linewidth]{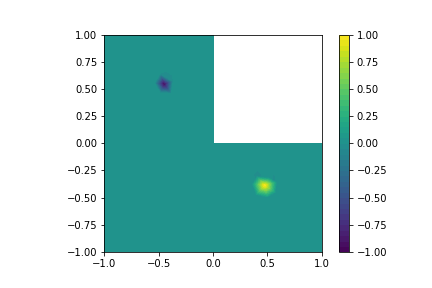}
        \caption{True sink and source.}
    \end{subfigure}
    \begin{subfigure}[b]{0.49\linewidth}        
        \centering
        \includegraphics[width=\linewidth]{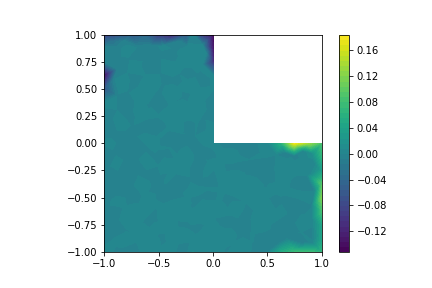}
        \caption{Inverse solution, $\alpha = 10^{-3}$.}
    \end{subfigure}\par
    \caption{\rem{Standard sparsity regularization (b) fails to recover the true sink-source configuration (a).}}
    \label{fig:motivation}
\end{figure}

\rem{
Figure \ref{fig:motivation}(b) shows the outcome of attempting to use boundary data and standard sparsity regularization to recover the "true" sink and source displayed in figure \ref{fig:motivation}(a). More precisely, we solved the problem 
\begin{equation*}
    \min_{\xx\in\Rn} \left\{\frac{1}{2}\|\AAA\xx - \bb\|_2^2 + \alpha\|\xx\|_1\right\}, 
\end{equation*}
where $\AAA$ is the forward matrix associated with \eqref{eq:in1}-\eqref{eq:in3.5} and $b$ denotes the Euclidean "version" of the boundary data $d$ generated by the true source-sink configuration displayed in figure \ref{fig:motivation}(a). Further information about the discretization of \eqref{eq:in1}-\eqref{eq:in3.5} and the matrix $\AAA$ is presented in the next section. 
}

\rem{
We observe that the classical LASSO scheme fails to provide adequate results for this problem, see figure \ref{fig:motivation}. Similar observations were made for all tested values of the regularization parameter $\alpha > 0$. In fact, results presented in 
\cite[pages 158--161]{burger2013inverse} and \cite{Elv21,Elv21c} reveal that internal sources will always be "moved" towards the boundary when standard regularization techniques are employed. This is the motivation for the present investigation.  
}

\section{Discretization} \label{sec:discretization}
In the finite element method (FEM) one usually approximates $f$ with a function $f_h$ in some finite element space $F_h$. The basis functions of $F_h$ typically have small and local support, e.g., they are "hat-functions" or characteristic functions of the involved grid cells, making them a suitable choice for incorporating spatial sparsity in our inverse problem. 

On the other hand, due to the complimentary condition \eqref{eq:in3.6}, 
we must search for a solution to our inverse problem in the subspace 
\[
\FhD = \left\{f_h \in F_h \; | \; \int_\Omega f_h = 0 \right\}
\]
of $F_h$. This subspace is spanned by the functions $$\psi_j = \phi_j - \frac{1}{| \Omega |} \int_\Omega \phi_j,$$
where $\{\phi_j\}$ is the set of basis functions for $F_h$. The imposed restriction of zero integral reduces the degrees of freedom of $\FhD$ by one, compared with the dimension of $F_h$. This makes the family $\{\psi_j\}$ a so-called \textit{frame} for $\FhD$. For an introduction to the use of frames in inverse problems, see, e.g., \cite{Chaux_2007}. 

We end this section with the definition of the forward matrix $\AAA$  associated with \eqref{eq:in3.5}. That is, $\AAA \in \mathbb{R}^{m \times n}$ represents the following chain of operations
\begin{equation*}
    \xx \in \mathbb{R}^n \quad \longrightarrow \quad f_h = \sum x_j \psi_j \quad \longrightarrow \quad u_h \quad \longrightarrow \quad T u_h|_{\partial\Omega} \in \mathbb{R}^m,
\end{equation*}
where $u_h = u_h(f_h)$ denotes the FEM solution of \eqref{eq:in1}-\eqref{eq:in3} generated by a given $f_h$ and $T$ yields the nodal values of $u_h|_{\partial\Omega}$. We discretized both $f_h$ and $u_h$ with first order Lagrange elements. Typically, $n \gg m$ and $\AAA$ will have a large null space. 


\section{Analysis} \label{sec:analysis}
Consider the equation
\begin{equation}\label{eq:1}
    \AAA\xx = \bb,
\end{equation}
where $\AAA \in \Rmn$ has a nontrivial null space $\nullspace{\AAA}$. We can multiply \eqref{eq:1} with the pseudo-inverse $\AAd$ of $\AAA$ to obtain
\begin{equation}\label{eq:2}
 \AAd\AAA\xx = \AAd\bb,
\end{equation}
where $\AAd\AAA$ is the orthogonal projection
\begin{equation}\label{eq:3}
    \PP = \AAd\AAA : \Rn \rightarrow \NAT.
\end{equation}

In \cite{Elv21} we introduced the diagonal {\em weight/regularization} matrix $\WW \in \Rnn$ defined by 
\begin{equation}\label{eq:4}
    \WW \ee_i = \|\PP\ee_i\|_2\ee_i,
\end{equation}
where we assume that $w_i=\|\PP\ee_i\|_2 \neq 0$ for $i=1,\, 2, \, \ldots,n$, i.e., none of the standard basis vectors belong to the null space of $\AAA$. 
Some of the beneficial mathematical properties of $\WW$ are explored in \cite{Elv21,Elv21b,Elv22,Elv21c}. In particular, let  $$\xx^\dagger = \AAA^\dagger \AAA \ee_j = \PP \ee_j$$ 
denote the minimum norm solution of $$\AAA \xx = \AAA \ee_j.$$ 
Then, 
 \begin{equation}\label{eq:max_property}
 \begin{split}
     &j = \argmax_{i \in \{1,2,...,n\}} |[\WW^{-1}\xx^\dagger]_i| = \argmax_{i \in \{1,2,...,n\}} |[\WW^{-1}\PP\ee_j]_i|, \\
     &[\WW^{-1}\PP\ee_j]_j = \|\PP\ee_j\|_2 \leq 1.
     \end{split}
 \end{equation}
 where, $[\mathbf{v}]_i$ represents the $i$'th component of the Euclidean vector $\mathbf{v}$. The property \eqref{eq:max_property} shows that $\WW^{-1} \xx^\dagger$ achieves its maximum for the "correct" index $j$, i.e., $\WW^{-1}\xx^\dagger$ can be used to recover the correct index $j$ from the "data" $\bb^\dagger = \AAA \ee_j$. 
 We will use this result at several occasions below. The proof of \eqref{eq:max_property} requires that the images under $\AAA$ of any two standard unit basis vectors are not parallel, i.e., 
\begin{equation} \label{NBP4.8}
    \AAA \ee_i \neq \eta \AAA \ee_j \quad \forall i, j \in \{ 1,2,\ldots,n \}, i\neq j, \forall \eta \in \mathbb{R}.
\end{equation}
 
 The property \eqref{eq:max_property} is the basic observation reported for the sLORETA algorithm in the EEG literature\footnote{In the sLORETA algorithm the weights are not defined in terms of the projection $\PP = \AAA^\dagger \AAA$, but instead one uses Tikhonov regularization, i.e., $(\AAA^T \AAA + \alpha \II)^{-1} \AAA^T \AAA \approx \AAA^\dagger \AAA$, to compute the weights. Hence, strictly speaking, \eqref{eq:max_property} may not hold for the sLORETA method unless $\alpha >0$ is very small.}  \cite{pascual2002standardized}. As far as the authors know, the first mathematical proof of \eqref{eq:max_property} was presented in connection with Theorem 4.1 in \cite{Elv21}, see also Theorem 2.1 in \cite{Elv21b}. 
 
\subsection{Basis pursuit} \label{sec:basis_pursuit}
This section is devoted to a study of the basis pursuit problem associated with \eqref{eq:1}, using the weight matrix \eqref{eq:4}:
\begin{equation} \label{NBP0}
  \min_\xx \| \WW \xx \|_1 \quad \textnormal{subject to} \quad \AAA \xx = \bb.  
\end{equation}
We want to analyze whether we can recover a sparse solution from exact data. That is, if the true data reads $$\bb = \bb^\dagger = \AAA\left(\sJee\right),$$ where $\mathcal{J} = \textnormal{supp}(\xx^*)$, can we recover
\begin{equation}\label{eq:5.1} 
\xx^* = \sJee 
\end{equation} 
by solving \eqref{NBP0}? 

\rem{
We first present, in Theorem \ref{thm_main_basis_pursuit} below, requirements which guarantee that: 
\begin{description}
\item[(a)] The true source-sink configuration $\xx^*$ is a solution of the basis pursuit problem \eqref{NBP0}. 
\item[(b)] The support of any other solution of this problem is contained in the support of $\xx^*$. 
\end{description}
}

\rem{
Part (a) is closely related to analogous results available in the classical sparsity regularization/LASSO literature, see, e.g., \cite{duval2017sparse,Fuchs04,grasmair10}. Typically, requirements similar to \eqref{NBP1}-\eqref{NBP2} in Theorem \ref{thm_main_basis_pursuit} are derived from the standard optimality conditions for the Lagrangian associated with \eqref{NBP0}, and the vector $\cc$ is referred to as the dual certificate\footnote{\rem{Employing the change of variable $\zz = \WW \xx$ in \eqref{NBP0}, we obtain the problem 
\begin{equation*} 
  \min_\zz \| \zz \|_1 \quad \textnormal{subject to} \quad \AAA \WW^{-1} \zz = \bb,   
\end{equation*}
which is in the standard form used in the LASSO literature.}}.
}

\rem{The support matter is analyzed in \cite{duval2017sparse} for the regularized (LASSO) counterpart to \eqref{NBP0}. More specifically, assuming that $\xx^*$ is the {\em only} solution of \eqref{NBP0}, the authors of \cite{duval2017sparse} study the support estimated by the LASSO approach when noise is added to the observation, which complements and extends the results presented in \cite{duval2015exact} and \cite{Fuchs04}. 
On the other hand, note that issue (b) concerns cases when \eqref{NBP0} potentially has {\em several} solutions. As far as the authors know, this situation has not received the same amount of attention: The proof of Theorem 2.2 in \cite{zhang2015necessary} shows that (b) holds, but property (b) is not explicitly mentioned/formulated in that paper. 
}

\rem{
Our analysis of (a) and (b) does not involve the Lagrangian or subgradients, and the argument for (b) differs significantly from the argument presented in \cite{zhang2015necessary}. Therefore, for the sake of completeness, we present a proof of Theorem \ref{thm_main_basis_pursuit} in appendix \ref{app:proof}.   
}
\begin{theorem}[The support will not increase and $\xx^*$ is a solution] \label{thm_main_basis_pursuit}
Assume that there exists a vector $\cc$ such that 
\begin{align}
\label{NBP1}
\frac{\PP \ee_i}{\| \PP \ee_i \|_2} \cdot \cc &= \textnormal{sgn}(x_i^*) \quad \forall i \in \JJ, \\
\label{NBP2}
\left| \frac{\PP \ee_i}{\| \PP \ee_i \|_2} \cdot \cc \right| &< 1 \quad \forall i \in \JJ^c. 
\end{align}
where $\JJ = \mathrm{supp} (\xx^*)$ and $\JJ^c=\{ 1,2,\ldots,n \} \setminus \JJ$.
Then any solution $\yy$ of 
\begin{equation}
\label{NBP3}
      \min_\xx \| \WW \xx \|_1 \quad \textnormal{subject to} \quad \AAA \xx = \AAA \overbrace{\sum_{\JJ} x_j^* \ee_j}^{= \xx^*}
\end{equation}
satisfies 
\[
\mathrm{supp} (\yy) \subseteq \mathrm{supp} (\xx^*), 
\]
and $\xx^*$ is a solution of \eqref{NBP3}. 
\end{theorem}

\rem{
We will now use Theorem \ref{thm_main_basis_pursuit} to prove that perfect recovery is possible in some particular cases for the discrete counterpart to the source-sink identification problem \eqref{eq:in3.5} discussed in the Introduction. Recall the definition \eqref{eq:4} of the diagonal weight matrix $\WW$, which involves the projection $\PP$ introduced in  equation \eqref{eq:3}. As mentioned above, some of the beneficial mathematical properties of $\WW$ for the single source case, using the screened Poisson equation as state equation, are revealed in \cite{Elv21,Elv21c}. 
}

\rem{
If $\AAA$ is the forward matrix associated with the model problem \eqref{eq:in1}-\eqref{eq:in3.5}, which involves the boundary value problem \eqref{eq:in1}-\eqref{eq:in3}, then the projection $\PP \ee_j = \AAA^\dagger \AAA \ee_j$ will typically not have local support: $\PP \ee_j$ is the minimum norm least squares solution of $\AAA \xx = \AAA \ee_j$ and can be approximated by applying a small amount of standard Tikhonov regularization. Moreover, a detailed analysis of the continuous counterpart to $\PP \ee_j$ is presented in \cite[section 2]{Elv21}, see also \cite[pages 158--161]{burger2013inverse}, and it turns out that the continuous counterpart to $\PP \ee_j$ achieves its maximum at the boundary $\partial \Omega$ and that most of its "mass/significant support" also is located close to $\partial \Omega$. 
}
\rem{
Consequently, the projections $\{ \PP \ee_j\}$ will typically be "almost disjoint" if the distances between the sources and the sinks are sufficiently large or if they are located in different "pockets" of the domain $\Omega$. Motivated by this discussion, we now prove that well-separated sources and sinks can be recovered, provided that the images under $\AAA$ of any two basis functions \rem{are not} parallel:
}
\begin{theorem}[Well-separated sources and sinks, i.e., projections with disjoint supports] \label{Co_disjoint_supports_BP}
Assume that \eqref{NBP4.8} holds. 
If 
\begin{equation} \label{NBP4.9}
\textnormal{supp}(\PP\ee_j) \cap \textnormal{supp}(\PP\ee_k) = \emptyset \mbox{ for all } j, k \in \JJ, j\neq k, 
\end{equation}
then
\[
\xx^* = \sJee, 
\] 
where $\JJ=\textnormal{supp}(\xx^*)$, is the unique solution of 
\begin{equation}
\label{NBP5}
      \min_{\xx \in \mathbb{R}^n} \| \WW \xx \|_1 \quad \textnormal{subject to} \quad \AAA \xx = \AAA \xx^*. 
\end{equation}
\end{theorem}
\begin{proof}
We will first show that \eqref{NBP1} and \eqref{NBP2} hold for an appropriate choice of $\cc$. Thereafter we prove the uniqueness. 

Let 
\begin{equation*}
    \cc = \sum_\JJ \textnormal{sgn}(x_j^*) \frac{\PP \ee_j}{\| \PP \ee_j \|_2} 
\end{equation*}
and note that \eqref{NBP4.9} assures that the projections $\{ \PP \ee_j \}_\JJ$ are orthogonal. 
Therefore, if $i \in \JJ$, then 
\begin{equation}
   \frac{\PP \ee_i}{\| \PP \ee_i \|_2} \cdot \cc = \frac{\PP \ee_i}{\| \PP \ee_i \|_2} \cdot \frac{\PP \ee_i}{\| \PP \ee_i \|_2} \textnormal{sgn}(x_i^*) = \textnormal{sgn}(x_i^*),  \label{NBP5.1}
\end{equation}
and we conclude that \eqref{NBP1} is fulfilled. 

Assume that $i \in \JJ^c$. Due to \eqref{NBP4.9} there can at most be one $k \in \JJ$ such that $i \in \textnormal{supp}(\PP \ee_{k})$. Consequently, provided that such a $k$ exists, and using the fact that $\PP$ is an \rem{{\em orthogonal projection}},   
\begin{align}
\nonumber
    \frac{\PP \ee_i}{\| \PP \ee_i \|_2} \cdot \cc &= \sum_\JJ \frac{\textnormal{sgn}(x_j^*)}{\| \PP \ee_i \|_2 \| \PP \ee_j \|_2} \, (\PP \ee_i,\PP \ee_j) \\
    \nonumber
    &=\sum_\JJ \frac{\textnormal{sgn}(x_j^*)}{\| \PP \ee_i \|_2 \| \PP \ee_j \|_2} \, (\ee_i,\PP \ee_j) \\ 
    \label{ny_likn1}
    &=\sum_\JJ \frac{\textnormal{sgn}(x_j^*)}{\| \PP \ee_i \|_2 \| \PP \ee_j \|_2} \, [\PP \ee_j]_i \\
    \nonumber
    &=\frac{\textnormal{sgn}(x_{k}^*)}{\| \PP \ee_i \|_2 \| \PP \ee_{k} \|_2} \, [\PP \ee_{k}]_i \\
    \nonumber
    &=\frac{\textnormal{sgn}(x_{k}^*)}{\| \PP \ee_i \|_2 \| \PP \ee_{k} \|_2} \, (\ee_i,\PP \ee_{k}) \\
    \nonumber
    &=\frac{\textnormal{sgn}(x_{k}^*)}{\| \PP \ee_i \|_2 \| \PP \ee_{k} \|_2} \, (\PP \ee_i,\PP \ee_{k}), 
\end{align}
where $[\PP \ee_j]_i$ denotes the $i$'th component of the (Euclidean) vector $\PP \ee_j$. 
Since \eqref{NBP4.8} holds, it follows from the Cauchy-Schwarz inequality and the fact that $i \in \JJ^c$ and $k \in \JJ$, i.e., $i \neq k$, that 
\begin{equation*}
    \left| \frac{\PP \ee_i}{\| \PP \ee_i \|_2} \cdot \cc \right| < 1,  
\end{equation*}
which is \eqref{NBP2}. 
On the other hand, if $i \in \JJ^c$ and $i \notin \textnormal{supp}(\PP \ee_j)$ for all $j \in \JJ$, then 
\begin{equation*}
    [\PP \ee_j]_i = 0 \quad \forall j \in \JJ. 
\end{equation*}
We thus instead find, see \eqref{ny_likn1}, that $$\frac{\PP \ee_i}{\| \PP \ee_i \|_2} \cdot \cc = 0,$$
and \eqref{NBP2} also holds in this case. 

Since both \eqref{NBP1} and \eqref{NBP2} are satisfied, Theorem \ref{thm_main_basis_pursuit} assures that $\xx^*$ solves \eqref{NBP5} and that any other solution $\yy$ of this problem must be such that $\textnormal{supp}(\yy) \subseteq \textnormal{supp}(\xx^*) = \JJ$. Hence, 
\begin{equation*}
    \AAA \yy = \AAA \xx^*
\end{equation*}
can be written in the form 
\begin{equation*}
    \AAA \sum_\JJ y_j \ee_j = \AAA \sum_\JJ  x_j^* \ee_j,
\end{equation*}
or, multiplying with $\AAA^\dagger$, 
\begin{equation*}
    \sum_\JJ y_j \PP \ee_j = \sum_\JJ  x_j^* \PP \ee_j. 
\end{equation*}
Finally, the orthogonality of the projections $\{ \PP \ee_j \}_\JJ$ implies that $y_j = x_j^*$ for all $j \in \JJ$, and it follows that $\xx^*$ is the unique solution of \eqref{NBP5}. 
\end{proof}

\rem{In this argument we used the fact that $\PP$ is an orthogonal projection. It is thus, as far as the authors know, an open problem whether Theorem \ref{Co_disjoint_supports_BP} also holds in a more general basis pursuit setting.}

Note that this theorem also covers the single-source-case, i.e., the case $\JJ = \{ j \}$. Then \eqref{NBP4.9} automatically holds, and the true source $\xx^* = x_j^* \ee_j$ can always be recovered. An alternative proof addressing the single-source-situation is presented in \cite{Elv21c}.  

\rem{
Theorem \ref{Co_disjoint_supports_BP}, and the discussion preceding it, suggest that non-convex domains are preferable to convex regions. We will return to this issue in the numerical examples section. Let us also mention that a rigorous analysis of the "almost disjoint projections case" is an open problem.
}

For the PDE constrained optimization problem \rem{\eqref{eq:in1}-\eqref{eq:in3.5}} it is plausible that a large number of sources or sinks located at the boundary, where data is recorded, can be recovered. This corresponds to the situation where the sources and sinks are close to, or in, the orthogonal complement of the null space of the associated forward operator/matrix. We will now not only prove that our methodology can handle such cases, but that also an additional interior source/sink can be detected.  

\begin{corollary}[Several sources and sinks in $\nullspace{\AAA}^\perp$, plus one more]
\label{Co_boundary_plusOne}
Let $\JJ$ be an index subset of $\{1, \,2, \, \ldots, \, n\}$ and assume that $\ee_{j} \in \nullspace{\AAA}^\perp$ for all $j \in \JJ$.
If \eqref{NBP4.8} holds, then 
\[
\xx^* = \sJee + x_{\notj}^* \ee_{\notj}
\] 
is the unique solution of 
\begin{equation}
\nonumber
      \min_{\xx \in \mathbb{R}^n} \| \WW \xx \|_1 \quad \textnormal{subject to} \quad \AAA \xx = \AAA \xx^*. 
\end{equation}
Here, $\notj$ is an index outside $\JJ$, i.e., $\notj \in \JJ^c$. 
\end{corollary}
\begin{proof}
First, 
\begin{equation*}
   \ee_{j} \in \nullspace{\AAA}^\perp \; \; \forall j \in \JJ \quad \Longrightarrow \quad  \PP \ee_j = \ee_j \; \;\forall j \in \JJ,
\end{equation*}
and it follows that 
\begin{equation} \label{NBP6}
    \textnormal{supp}(\PP \ee_j) \cap \textnormal{supp}(\PP \ee_k) = \emptyset \quad \forall j,k \in \JJ, j \neq k.
\end{equation}

Second, for any $j \in \JJ$, keeping in mind that $\PP^T = \PP$, 
\begin{equation*}
    [\PP \ee_\notj]_j = (\PP \ee_\notj, \ee_j) = (\ee_\notj, \PP \ee_j) = (\ee_\notj, \ee_j) = 0. 
\end{equation*}
Consequently,  
\begin{equation} \label{NBP7}
    \emptyset = \textnormal{supp}(\PP \ee_\notj) \cap \textnormal{supp}(\ee_j) = \textnormal{supp}(\PP \ee_\notj) \cap \textnormal{supp}(\PP \ee_j) \quad \forall j \in \JJ. 
\end{equation}

Finally, \eqref{NBP6} and \eqref{NBP7} imply that 
\begin{equation*} 
    \textnormal{supp}(\PP \ee_j) \cap \textnormal{supp}(\PP \ee_k) = \emptyset \quad \forall j,k \in \JJ \cup \{ \notj \}, j \neq k, 
\end{equation*}
and the result therefore follows from Theorem \ref{Co_disjoint_supports_BP}.
\end{proof}

\rem{As for Theorem \ref{Co_disjoint_supports_BP}, we do not know whether Corollary \ref{Co_boundary_plusOne} holds in a more general basis pursuit setting: The proof above explicitly employs that $\PP$ is an orthogonal projection and invokes Theorem \ref{Co_disjoint_supports_BP} which also assumes that $\PP$ is a projection.}

In practice it seems important that a source is not misinterpreted as a sink or vice versa. The next result addresses this issue. 
\begin{theorem}[\rem{A source will not misinterpreted as a sink or vice versa}] 
\rem{
Assume that there exists a vector $\cc$ such that \eqref{NBP1} and \eqref{NBP2} hold with $\JJ = \textnormal{supp}(\xx^*)$. 
Then any solution $\yy$ of 
\begin{equation} \label{sameSign2}
    \min_{\xx} \| \WW \xx \|_1 \quad \textnormal{subject to} \quad \AAA \xx = \AAA\xx^*
\end{equation}
obeys   
\begin{equation*}
    \textnormal{sgn} (y_k) \in \{0, \textnormal{sgn} (x^*_k)\} \quad \forall k \in \{1,2, \ldots, n\}.
\end{equation*}
}
\end{theorem}
\rem{
\begin{proof}
First, note that \eqref{NBP1}, because $\PP^T=\PP$, can be written in the form
\begin{equation}\label{NBP1.rev}
    [\PP\cc]_i = \|\PP\ee_i\|_2\textnormal{sgn}(x_i^*), \quad \forall i \in \JJ.
\end{equation}
Next, define $\qq = \xx^* - \yy \in \nullspace{\AAA}$. From Theorem \ref{thm_main_basis_pursuit} we have that $\textnormal{supp}(\yy) \subseteq \textnormal{supp}(\xx^*)$, and consequently, $\textnormal{supp}(\qq) \subseteq \textnormal{supp}(\xx^*)$. Combining this observation with \eqref{NBP1.rev} and the fact that $\xx^*$ solves \eqref{sameSign2}, we obtain
\begin{eqnarray*}
    \PP\cc \cdot \qq &=& \sum_\JJ \|\PP\ee_i\|_2\textnormal{sgn}(x_i^*)q_i \\
    &=& \sum_\JJ \|\PP\ee_i\|_2\textnormal{sgn}(x_i^*)(x_i^* - y_i) \\
    &=& \|\WW\xx^*\|_1 - \sum_\JJ \|\PP\ee_i\|_2\textnormal{sgn}(x_i^*)y_i \\
    &\geq& \|\WW\xx^*\|_1 - \|\WW\yy\|_1 = 0,
\end{eqnarray*}
with equality only if $\textnormal{sgn}(y_i) \in \{0,\textnormal{sgn}(x_i^*)\}$ for all $i \in \JJ$. The estimate must indeed hold with equality since $\qq \in \nullspace{\AAA} = \nullspace{\PP}$ and thus
$\PP\cc\cdot \qq = \cc \cdot \PP\qq = 0$. This completes the proof.
\end{proof}
}
\rem{We next prove that $\xx^*$ is the only solution to the basis pursuit problem if and only if $\AAA$ is injective on $\textnormal{supp} (\xx^*)$. This result is not new, see \cite{Fuchs04,grasmair10, zhang2015necessary}. Nevertheless, since we know that the support of any solution $\yy$ of \eqref{NBP3} must be a subset of the support of $\xx^*$, see Theorem \ref{thm_main_basis_pursuit}, our proof of this fact becomes short.}
Note that we use the notation 
\begin{equation*}
    \mathbf{sgn}(\vv) = \left( \textnormal{sgn}(v_1), \textnormal{sgn}(v_2), \ldots, \textnormal{sgn}(v_n) \right)^T
\end{equation*}
\rem{in the proof below.}

\rem{
\begin{theorem}[Uniqueness] \label{BP_uniqueness_new} 
Assume that $\xx^*$ is such that \eqref{NBP1} and \eqref{NBP2} hold. Then $\xx^*$ is the unique solution of \eqref{NBP3} if and only if $\AAA$ is injective on the support of $\xx^*$. 
\end{theorem}
}
\begin{proof}
\rem{
According to Theorem \ref{thm_main_basis_pursuit}, $\xx^*$ solves \eqref{NBP3} and any other solution $\yy$ of \eqref{NBP3} satisfies $\textnormal{supp} (\yy) \subseteq \textnormal{supp} (\xx^*)$. The constraint in \eqref{NBP3} yields that $\AAA \yy = \AAA \xx^*$, and consequently, if $\AAA$ is injective on the support of $\xx^*$, \eqref{NBP3} can not have other solutions than $\xx^*$. 
}

\rem{
Assume that $\xx^*$ is the only solution to \eqref{NBP3} and that $\AAA$ is not injective on $\textnormal{supp} (\xx^*)$. Then there exists a unit vector $\qq \in \nullspace{\AAA}$ with $\textnormal{supp} (\qq) \subseteq \textnormal{supp} (\xx^*)$. The directional subderivative $\partial_\qq \| \WW \xx \|_1$ at $\xx^*$ must obey, because $\xx^*$ is a minimizer,  
\begin{eqnarray*}
    \partial_\qq \| \WW \xx^* \|_1 = \mathbf{sgn}(\xx^*)^T \WW \qq = 0. 
\end{eqnarray*}
(Actually, $\| \WW \xx \|_1$ has a standard directional derivative at $\xx^*$ in the direction $\qq$). 
Since $\textnormal{supp} (\qq) \subseteq \textnormal{supp} (\xx^*)$, there exists $\hat{t} > 0$ such that 
\begin{eqnarray*}
\mathbf{sgn}(\xx^* + t \qq)=\mathbf{sgn}(\xx^*) \quad \forall t \in [0,\hat{t}].    
\end{eqnarray*}
This implies that 
\begin{eqnarray*}
    \partial_\qq \| \WW (\xx^*+t \qq) \|_1 = \mathbf{sgn}(\xx^* + t \qq)^T \WW \qq = \mathbf{sgn}(\xx^*)^T \WW \qq = 0 \quad \forall t \in [0,\hat{t}],
\end{eqnarray*}
and we conclude that $\xx^*+t \qq$, for any $t \in [0,\hat{t}]$, also would solve \eqref{NBP3}. Thus, if $\xx^*$ is the unique solution of \eqref{NBP3}, $\AAA$ must be injective on $\textnormal{supp} (\xx^*)$.
} 
\end{proof}

\rem{Using modern software tools, it is "easy" to check, for a given true sink-source vector $\xx^*$, whether $\AAA$ is injective on $\textnormal{supp}(\xx^*)$.} 



\subsection{Regularized problems} \label{sec:regularized_problems}
In this section we establish regularized counterparts to Theorem \ref{Co_disjoint_supports_BP} and Corollary \ref{Co_boundary_plusOne}. This is important because in practice one needs to apply regularization in order to avoid disastrous amplification of noise in real world data.  

If we apply weighted regularization to \eqref{eq:2}, using the weight matrix $\WW$, we obtain
\begin{equation}\label{eq:5}
    \min_{\xx\in\Rn} \left\{\frac{1}{2}\|\AAd\AAA\xx - \AAd\bb\|_2^2 + \alpha\|\WW\xx\|_1\right\}.
\end{equation}
Since $\WW$ is defined in terms of the projection $\PP = \AAd\AAA$, it turns our that it is convenient to consider \eqref{eq:5} instead of the standard formulation which is based on \eqref{eq:1}. In fact, \rem{except for Theorem \ref{thm:convergence}}, we have not succeeded in proving similar results, to those presented in this subsection, for the problem 
\begin{equation*}
    \min_{\xx\in\Rn} \left\{\frac{1}{2}\|\AAA\xx - \bb\|_2^2 + \alpha\|\WW\xx\|_1\right\}.
\end{equation*}

If the underlying problem is ill posed, then $\AAA$ will typically have very small non-zero singular values, and the use of $\AAA^\dagger$ is not recommendable. In practice we therefore replace $\AAA^\dagger$ with a more "well behaved" matrix. This can be accomplished by applying truncated SVD or standard Tikhonov regularization. 

Analogously to the investigation of the basis pursuit problem, we will analyze whether $$\xx^* = \sJee$$ can be approximately recovered from the true data $$\bb=\bb^\dagger = \AAA\left(\sJee\right).$$ 
For this problem, using the fact that $\PP = \AAd\AAA$, the minimization problem \eqref{eq:5} becomes
\begin{equation}\label{eq:6}
    \yy^*_\alpha = \argmin_{\xx\in\Rn} \left\{\frac{1}{2}\left\|\PP\xx-\PP\left(\sJee\right)\right\|_2^2 + \alpha\|\WW\xx\|_1\right\}.
\end{equation}

\rem{
Lemma \ref{lem:a_implies_c} and Theorem \ref{thm:main} below yield two criteria for determining the {\em approximate} identifiability of $\xx^*$. (Approximate in the sense that the support of $\xx^*$ is preserved in the regularized solution and that one obtains the "correct solution" as $\alpha \rightarrow 0$). More specifically, it turns out that it is sufficient to explore whether the $s \times s$ linear system \eqref{eq:7.1} below, where $s$ is the cardinality of $\JJ = \textnormal{supp}(\xx^*)$, admits a solution and whether the product of the involved system matrix and this solution has components, associated with indexes in $\JJ^c = \{1,\,2,\, \ldots, \, n \} \setminus \JJ$, belonging to the interval $(-1,1)$, see \eqref{eq:7.2}. These two criteria are stronger that \eqref{NBP1} and \eqref{NBP2}, which we now prove. Recall that we use the notation $[\vv]_i$ for the $i$'th component of the Euclidean vector $\vv$. 
}
\rem{
\begin{lemma} \label{lem:a_implies_c}
    Let $\JJ = \textnormal{supp}(\xx^*)$, with cardinality $s$, and assume that there exist $a_1, a_2, \ldots, a_s$ such that  
    \begin{align}
        \left[\sJ \aam \WW^{-1}\PP\ee_j\right]_i &= \textnormal{sgn} (x_i^*) \quad \forall i \in \JJ 
        \label{eq:7.1}\tag{C.1} 
        \end{align}
        and 
        \begin{align}
        \left[\sJ \aam \WW^{-1}\PP\ee_j\right]_i &\in \rem{(-1, 1)} \quad \forall i \in \JJ^c.  \label{eq:7.2}\tag{C.2}
    \end{align}
    Then there exists a vector $\cc$ such that \eqref{NBP1} and \eqref{NBP2} hold.
\end{lemma}
}
\rem{
\begin{proof}
In the following chain of equalities we use the facts that $\PP^T=\PP$ ($\PP$ is an orthogonal projection), $\WW^{-1} \ee_i = \| \PP \ee_i \|_2^{-1} \ee_i$ (see \eqref{eq:4}) and that $\WW^{-1}$ is a diagonal matrix 
 \begin{align*}
    \frac{\PP \ee_i}{\| \PP \ee_i \|_2} \cdot \cc &= \frac{\ee_i}{\| \PP\ee_i \|_2} \cdot \PP \cc \\
    &= \WW^{-1} \ee_i \cdot \PP \cc \\
    &= \ee_i \cdot \WW^{-1} \PP \cc \\
    &= \ee_i \cdot \sum_{\JJ \cup \JJ^c} c_j \WW^{-1} \PP \ee_j \\
    &= \left[ \sum_{\JJ \cup \JJ^c} c_j \WW^{-1} \PP \ee_j \right]_i,
 \end{align*}   
 where we have used the notation $\cc=\sum_{\JJ \cup \JJ^c} c_j \ee_j$. Consequently, if \eqref{eq:7.1} and \eqref{eq:7.2} hold, then by choosing 
 \begin{equation*}
     c_j=\left\{ \begin{array}{cc}
         a_j & j \in \JJ \\
         0 &  j \in \JJ^c
     \end{array} \right.,
 \end{equation*}
 we find that \eqref{NBP1} and \eqref{NBP2} are satisfied.
\end{proof}
}

\begin{theorem}[\rem{Necessary and } sufficient criteria] \label{thm:main}
    Let $\PP$ and $\WW$ be the matrices defined in \eqref{eq:3} and \eqref{eq:4}, respectively, and let $\JJ = \textnormal{supp}(\xx^*)$. If the linear system \rem{\eqref{eq:7.1}
        has a solution $(a_1, a_2, \ldots, a_s)$ which satisfies \eqref{eq:7.2},}
    then $$\yy^*_\alpha = \sum_{\JJ}(x_j^* - \alpha \aam) \ee_j$$ solves \eqref{eq:6}, provided that $\alpha > 0$ is chosen such that  
    \begin{equation}
    \label{eq:7.25}
        \textnormal{sgn}(x_j^* - \alpha \aam) = \textnormal{sgn}(x_j^*) \quad \forall j \in \JJ. 
    \end{equation}
    \rem{Furthermore, the solution $\yy^*_\alpha$ of \eqref{eq:6} is \textit{unique} if and only if $\AAA$ is injective on the support of $\xx^*$.}
\end{theorem}
\begin{remark}
Before proving the theorem, we note that the conditions \eqref{eq:7.1} and \eqref{eq:7.2} only involve the sign of $x_j^*, \, j \in \JJ$. This leads to a rather 
interesting observation: Whether we are, in principle, able to \rem{\em approximately} recover an $s$-sparse solution is independent of the magnitudes $\{ | x_j^* | \}_\JJ$ of the individual sources and sinks. Moreover, when $\alpha > 0$ is sufficiently small, \eqref{eq:7.25} will always hold because we consider finite dimensional problems. We also mention that the linear system \eqref{eq:7.1} typically will be small, unless the number $s$ of sparse sinks and sources is large. In the corollaries below we explore situations where we can guarantee that \textit{both} \eqref{eq:7.1} and \eqref{eq:7.2} hold. 

This theorem provides us with a simple two step strategy for checking whether a configuration of sources and sinks can be \rem{\em approximately} recovered under ideal conditions: 
\begin{itemize}
    \item Solve, if possible, \eqref{eq:7.1} for a given configuration of sinks and sources.
    \item If \eqref{eq:7.1} is solvable, verify whether \eqref{eq:7.2} holds. 
\end{itemize}
For example, 
one can investigate whether certain configurations of sinks and sources 
are detectable, provided that the noise level is sufficiently small. 
\end{remark}
\begin{proof}[Proof of Theorem \ref{thm:main}] $ $ \newline 
Let us define the cost-functional $\mathcal{T}_\alpha: \Rn \rightarrow \mathbb{R}$ associated with \eqref{eq:6},   \begin{equation} \label{eq:8.1}
      \mathcal{T}_\alpha(\xx) = \underbrace{\frac{1}{2}\left\|\PP\xx-\PP\left(\sJee\right)\right\|_2^2}_{=g(\xx)}
      + \underbrace{\alpha\|\WW\xx\|_1}_{=\alpha h(\WW\xx)},
  \end{equation}
  where $g(\cdot)$ and $h(\WW\cdot)$ represent the fidelity and regularization terms, respectively.
  According to standard convex optimization theory, $\xx$ is a minimizer of $\mathcal{T}_\alpha$ if and only if
  \begin{eqnarray*}
      \mathbf{0} &\in& \partial\mathcal{T}_\alpha(\xx) \\
      &=& \nabla g(\xx) + \alpha\WW^T\partial h(\WW\xx), 
  \end{eqnarray*}
  where "$\partial$" denotes the subgradient. 
  Since $\WW^T = \WW$, we can multiply with $\WW^{-1}$ to obtain
  \begin{equation*}
     -\WW^{-1}\nabla g(\xx) \in \alpha\partial h(\WW\xx),
  \end{equation*}
  and from the expression \eqref{eq:8.1} for $g$ we find, keeping in mind that $\PP^T \PP = \PP \PP = \PP$, 
  \begin{equation}\label{eq:8.2}
     \WW^{-1}\PP\left(\sJee - \xx\right) \in \alpha \partial h(\WW\xx).
  \end{equation}  
  We also observe, using the fact that $h(\zz)=\| \zz \|_1$ and that $\WW$ is a diagonal matrix with positive entries at its diagonal,  
  \begin{equation*}
     [\partial h(\WW\xx)]_i =  [\partial h(\WW [x_1 \, x_2 \, \ldots \, x_n]^T)]_i =  
        \begin{cases} 
            \{1\}, & x_i > 0, \\
            \{-1\}, & x_i < 0, \\
            [-1,1], & x_i = 0.
        \end{cases}
  \end{equation*}
  
  We will now investigate whether there exist scalars $\{\tilde{\gamma}_j\}_\mathcal{J}$ such that $\xx = \sJgtee$, $\gamt \neq 0$, satisfies the optimality criterion  \eqref{eq:8.2}. Note that 
  \begin{equation} \label{eq:8.3}
  \left[\partial h\left(\WW\sJgtee\right)\right]_i =         
        \begin{cases} 
            \{ \textnormal{sgn}(\tilde{\gamma}_j) \}, & i \in \mathcal{J}, \\
            [-1,1], & i \notin \mathcal{J},
        \end{cases}
  \end{equation}
  so the condition \eqref{eq:8.2}, with $\xx = \sJ\gamt\ee_j$, becomes   
  \begin{equation} \label{eq:8.35}
      \left[\sJ (x_j^*-\gamt) \WW^{-1}\PP\ee_j\right]_i \in \alpha
        \begin{cases} 
            \{ \textnormal{sgn}(\tilde{\gamma}_j) \}, & i \in \mathcal{J}, \\
            [-1,1], & i \notin \mathcal{J}.
        \end{cases}
  \end{equation}
  Assume that \eqref{eq:7.1}-\eqref{eq:7.2} hold. Then, with $$\gamt = x_j^* - \alpha \aam,$$ we can conclude from \eqref{eq:7.1}-\eqref{eq:7.2} that
    \begin{equation} \label{eq:8.4}
      \left[\sJ (x_j^*-\gamt) \WW^{-1}\PP\ee_j\right]_i = 
      \alpha \left[\sJ \aam \WW^{-1}\PP\ee_j\right]_i \in \alpha
        \begin{cases} 
            \{ \textnormal{sgn}(x_j^*) \}, & i \in \JJ, \\
            \rem{(-1,1)}, & i \in \JJ^c,
        \end{cases}
    \end{equation}
    and \eqref{eq:8.35} is satisfied when $\alpha > 0$ is chosen such that $\textnormal{sgn}(\tilde{\gamma}_j) = \textnormal{sgn}(x_j^*)$, $\forall j \in \JJ$, i.e., when \eqref{eq:7.25} holds. \rem{Hence, we have proved that $\yy^*_\alpha = \sum_{\JJ}(x_j^* - \alpha \aam) \ee_j$ solves \eqref{eq:6}.}  

    \rem{
    According to Lemma \ref{lem:a_implies_c}, if there exist $a_1,a_2,\ldots,a_s$ such that \eqref{eq:7.1} and \eqref{eq:7.2} are satisfied, then there exists a dual certificate $\cc$ such that \eqref{NBP1} and \eqref{NBP2} hold. It therefore follows from \cite[Theorem 2.1]{zhang2015necessary} that $\yy^*_\alpha$ is the unique solution of \eqref{eq:6} if and only if $\AAA$ is injective on the support of $\xx^*$. (In order to invoke \cite[Theorem 2.1]{zhang2015necessary} one must use the change of variable $\zz=\WW \xx$). 
    }
\end{proof}  

The regularized counterpart to Theorem \ref{Co_disjoint_supports_BP}, i.e., the case with well-separated sinks and sources, reads: 
\begin{corollary}[Well-separated sources and sinks, i.e., disjoint projections] 
\label{reg_Th_disjoint_projections}
Let $\xx^* = \sJee$, $\JJ = \textnormal{supp}(\xx^*)$, and assume that \eqref{NBP4.8} holds. If $$\textnormal{supp}(\PP\ee_j) \cap \textnormal{supp}(\PP\ee_k) = \emptyset \quad j\neq k, j, \, k \in \mathcal{J},$$ then the conditions \eqref{eq:7.1}-\eqref{eq:7.2} hold. Furthermore, $$\yy^*_\alpha = \sum_{\JJ}(x_j^* - \alpha \aam) \ee_j$$ solves \eqref{eq:6}, provided that $\alpha > 0$ is chosen such that \eqref{eq:7.25} is satisfied.
\end{corollary}
\begin{proof}
 Since $\WW$ is a diagonal matrix, it follows that $\textnormal{supp}(\WW^{-1} \PP\ee_k) = \textnormal{supp}(\PP\ee_k)$ for all $k$. Consequently, 
 \begin{equation}\label{eq:10}
 \textnormal{supp}(\WW^{-1} \PP\ee_j) \cap \textnormal{supp}(\WW^{-1} \PP\ee_k) = \emptyset \quad j\neq k, \, j, k \in \mathcal{J}.
 \end{equation}
 As proven in \cite{Elv21b}, see \eqref{eq:max_property}, because \eqref{NBP4.8} is assumed to hold,  $\WW^{-1}\PP\ee_j$ achieves its unique absolute maximum for index $j$, i.e., $\left[\WW^{-1}\PP\ee_j\right]_j \neq 0$. Due to \eqref{eq:10}, condition \eqref{eq:7.1} therefore simplifies to 
 \begin{equation} \label{eq:11}
    \left[\aam \WW^{-1}\PP\ee_j\right]_j = \textnormal{sgn}(x_j^*) \quad \forall j \in \mathcal{J},
 \end{equation}
 which clearly can be solved for $\{ \aam \}_\JJ$.
 
 Furthermore, again because of the non-overlapping supports \eqref{eq:10}, there can for each $i \in \JJ^c$ be at most one element $\iota(i) \in \JJ$ such that $[\PP\ee_{\iota(i)}]_i \neq 0$. Consequently, condition \eqref{eq:7.2} reads
     \begin{equation} \label{eq:12}
        \left[a_{\iota(i)} \WW^{-1}\PP\ee_{\iota(i)}\right]_i \in (-1,1) \quad \forall i \in \JJ^c.
    \end{equation}
 Equation \eqref{eq:11} implies that $\left| \left[a_{\iota(i)} \WW^{-1}\PP\ee_{\iota(i)}\right]_{\iota(i)} \right| = 1$, and 
 the maximum property \eqref{eq:max_property} therefore implies that \eqref{eq:12} is satisfied.
\end{proof}

The next results shows that a regularized version of Corollary \ref{Co_boundary_plusOne} holds: 
\begin{corollary}[Several sources and sinks in $\nullspace{\AAA}^\perp$, plus one more]
\label{reg_Th_boundary_plusOne}
 Assume that \eqref{NBP4.8} holds. 
 If $\ee_j \in \NAT, \forall j \in \mathcal{J}\setminus\{\notj\}$, then the conditions \eqref{eq:7.1}-\eqref{eq:7.2} are satisfied, and $$\yy^*_\alpha = \sum_{\JJ}(x_j^* - \alpha \aam) \ee_j$$ solves \eqref{eq:6}, provided that $\alpha > 0$ is chosen such that \eqref{eq:7.25} is satisfied.
 \end{corollary}
 \begin{proof}
    Using Theorem \ref{thm:main}, we only need to show that \eqref{eq:7.1}-\eqref{eq:7.2} hold.
    First, note that for $\ee_j \in \NAT$ it follows from the definition \eqref{eq:3} of the orthogonal projection $\PP$ that $\PP\ee_j = \ee_j$. Therefore, $$\WW^{-1}\PP\ee_j = \WW^{-1} \ee_j = \ee_j,$$ where the last equality follows from \eqref{eq:4}: 
    \begin{equation*}
        \ee_j \in \NAT \Rightarrow \| \PP \ee_j \|_2 = \| \ee_j \|_2 = 1. 
    \end{equation*}
    Consequently, \eqref{eq:7.1} simplifies to
    \begin{equation}\label{eq:10.1}
      \left[\sum_{\mathcal{J}\setminus{\notj}} \aam \ee_j\right]_i + a_{\notj}[\WW^{-1}\PP\ee_{\notj}]_i = \textnormal{sgn}(x_i^*) \quad \forall i \in \JJ.
    \end{equation}
    This can be written as the matrix-vector equation
    \begin{equation*}
        \begin{bmatrix}
            1 & 0 & \cdots & 0 & [\WW^{-1}\PP\ee_\notj]_{j_1} & 0 & \cdots & 0 \\
            0 & 1 & \cdots & 0 & [\WW^{-1}\PP\ee_\notj]_{j_2} & 0 & \cdots & 0 \\
              &   &        &   &   \vdots                     &   &        &   \\
            0 & 0 & \cdots & 0 & [\WW^{-1}\PP\ee_\notj]_{\notj} & 0 & \cdots & 0 \\
              &   &        &   &   \vdots                     &   &        &   \\
            0 & 0 & \cdots & 0 & [\WW^{-1}\PP\ee_\notj]_{j_s} & 0 & \cdots & 1
        \end{bmatrix}
        \mathbf{a} =
    \begin{bmatrix}
        \textnormal{sgn}(x_{j_1}^*) \\
        \textnormal{sgn}(x_{j_2}^*) \\
        \vdots\\
        \textnormal{sgn}(x_{\notj}^*)\\
        \vdots \\
        \textnormal{sgn}(x_{j_s}^*)
    \end{bmatrix}
    \end{equation*}
    From the maximum property \eqref{eq:max_property} we know that 
    \[
    |[\WW^{-1}\PP\ee_{\notj}]_{j_q}| < |[\WW^{-1}\PP\ee_{\notj}]_{\notj}| \leq 1, \quad j_q \neq \notj. 
    \]
    Consequently, the matrix above is strictly diagonally dominant, and this system therefore has a unique solution. Hence, condition \eqref{eq:7.1} is satisfied.
    
    To show that \eqref{eq:7.2} also holds, we first note that for $i = \notj$, \eqref{eq:10.1} asserts that $a_\notj$ must satisfy
    \begin{equation*}
        a_\notj[\WW^{-1}\PP\ee_{\notj}]_{\notj} = \textnormal{sgn}(x_\notj^*).
    \end{equation*}
    Again, by \eqref{eq:max_property} we then have that for $i \in \JJ^c$,
    \begin{equation*}
        a_\notj[\WW^{-1}\PP\ee_{\notj}]_{i} \in (-1,1).
    \end{equation*}    
    In this case, this is sufficient to conclude that condition \eqref{eq:7.2} holds because $\PP \ee_j = \ee_j$ for all $j \in \mathcal{J}\setminus\{\notj\}$.
\end{proof}

\rem{
We have so far derived explicit expressions for the inverse regularized solution in some special cases. We end this section with a more general convergence result for cases involving noise, which, by employing the norm
\begin{equation}
    \|\cdot\|_\WW = \|\WW\cdot\|_2,
\end{equation}
can be derived almost immediately from results presented in \cite{grasmair10}:
\begin{theorem}\label{thm:convergence}
    Let the true data $\bb^\dagger$ be generated by $\bb^\dagger = \AAA\xx^*$ for some $\xx^*$. Assume that $\xx^*$ satisfies \eqref{NBP1}-\eqref{NBP2} and that $\AAA$ is injective on the support of $\xx^*$. Denote
    \begin{equation}\label{eq:formA}
        \xx_\alpha = \argmin_{\xx\in\Rn} \left\{\frac{1}{2}\|\AAA\xx - \bb^\delta\|_2^2+\alpha\|\WW\xx\|_1\right\},
    \end{equation}
    where $\bb^\delta$ is a noisy measurement of $\bb^\dagger$ with the bound $\|\bb^\delta-\bb^\dagger\|_2 \leq \delta$. Then, for every $C > 0$ there exists $c > 0$ such that 
    \begin{equation}\label{eq:convX}
        \|\xx_\alpha - \xx^*\|_{\WW} \leq c\delta,
    \end{equation}
    for the choice $\alpha = C\delta$.
\end{theorem}
}
\rem{
\begin{proof}
    From Theorem \ref{BP_uniqueness_new}, we know that $\xx^*$ is the unique solution of 
    \begin{equation*}
        \min_\xx \|\WW\xx\|_1 \quad \textnormal{subject to} \quad \AAA\xx = \bb^\dagger
    \end{equation*}
    because we assume that $\AAA$ is injective on the support of $\xx^*$. 
    Consequently, since $\WW$ is non-singular, $\zz^* = \WW\xx^*$ is the unique solution of
    \begin{equation} \label{eq:unique_z}
        \min_\zz \|\zz\|_1 \quad \textnormal{subject to} \quad \AAA\WW^{-1}\zz = \bb^\dagger.
    \end{equation}
    The associated variational formulation with noisy data reads 
    \begin{equation*}
        \zz_\alpha = \argmin_\zz \left\{\frac{1}{2}\|\AAA\WW^{-1}\zz - \bb^\delta\|_2^2+\alpha\|\zz\|_1\right\}. 
    \end{equation*}
    Since the solution $\zz^*$ of \eqref{eq:unique_z} is unique, it follows from \cite[Lemma 4.5]{grasmair10} that \cite[Condition 4.3]{grasmair10} holds and thus \cite[Theorem 4.7]{grasmair10} implies that 
    \begin{equation}
        \|\zz_\alpha - \zz^*\|_2 \leq c\delta. \label{eq:convZ}
    \end{equation}
    Recall the definition \eqref{eq:formA} of $\xx_\alpha$.  Then, since $\zz_\alpha = \WW\xx_\alpha$ and $\zz^* = \WW\xx^*$, it follows from \eqref{eq:convZ} that 
    \begin{equation*}
        \|\WW\xx_\alpha - \WW\xx^*\| = \|\xx_\alpha - \xx^*\|_\WW \leq c\delta.
    \end{equation*}
\end{proof}
}

\rem{
Recall that the results presented prior to Theorem \ref{thm:convergence} in this subsection concern \eqref{eq:6}, with $\PP = \AAA^\dagger \AAA$, and not \eqref{eq:formA}. As mentioned earlier, since the underlying problem is ill posed, it is not advisable to apply $\AAA^\dagger$ in practise. We therefore now briefly consider the formulation
\begin{equation}\label{eq:formAd}
    \yy_\alpha = \argmin_{\xx\in\Rn} \left\{\frac{1}{2}\|\AAA_k^\dagger\AAA\xx - \AAA_k^\dagger\bb^\delta\|_2^2 + \alpha\|\WW\xx\|_1\right\},
\end{equation}
where $\AAA_k^\dagger$ represents the "SVD-truncated" approximation of $\AAA^\dagger$. ($k$ is the number of non-zero singular values included in the approximation). Clearly, 
\begin{equation*}
    \|\AAA_k^\dagger(\bb^\delta - \bb^\dagger)\|_2 \leq \|\AAA_k^\dagger\|\|\bb^\delta - \bb^\dagger\|_2 = \|\AAA_k^\dagger\|\delta,
\end{equation*}
and we therefore, in this case, get the stability estimate 
\begin{equation}\label{eq:convY}
    \|\yy_\alpha - \xx^*\|_{\WW} \leq c\|\AAA_k^\dagger\|\delta
\end{equation}
instead of \eqref{eq:convX}. 
This again follows from the theory developed in \cite{grasmair10}, provided that \eqref{NBP1} and \eqref{NBP2} hold for $\PP_k=\AAA_k^\dagger\AAA$ and that $\PP_k$ is injective on the support of $\xx^*$. 
Notice that, without a truncation of $\AAA^\dagger$, \eqref{eq:convY} becomes meaningless in the infinite dimensional setting because then the pseudo-inverse is unbounded. How to determine the optimal size of the truncation parameter $k$ and whether \eqref{eq:formA} or \eqref{eq:formAd} are preferable for computing  inverse solutions, for a given noise level $\delta > 0$, are open problems.
}


\section{Numerical examples} \label{sec:numerical_experiments}
In this section we illuminate the theoretical results presented in section \ref{sec:regularized_problems}. More specifically, we will show experiments where the assumptions needed in corollaries \ref{reg_Th_disjoint_projections} and \ref{reg_Th_boundary_plusOne} are {\em almost} satisfied in order to demonstrate the robustness of the method. (We do not present computations for problems satisfying the assumptions exactly because then we know that the true sinks and sources will be perfectly recovered.)   

Concerning the numerical solution of \eqref{eq:6}, where we recall that $\PP=\AAA^\dagger \AAA$, we employed truncated SVD to obtain a "well-behaved" approximation of $\AAA^\dagger$, \rem{using $50$ non-zero singular values in examples $0$ and $1$, $20$ non-zero singular values in example $2$ and $10$ non-zero singular values in example 3}. 
\rem{The conductivity $\sigma$ equaled $1$ in examples $0$ and $1$, cf. equation \eqref{eq:in1}, whereas in the examples $2$ and $3$ we used $\sigma(x,y) = 2 + \sin(x)\cos(y)$}. All domains were partitioned in terms of non-uniform grids, and the mesh parameter $h$ in the inverse computations was twice as large as in the forward simulations.  All code was written in MATLAB and Python, using the FEniCS and Scipy libraries for the latter.

\subsection{\rem{Example 0}} 
\begin{figure}[H]
        \centering
        \includegraphics[width=0.49\linewidth]{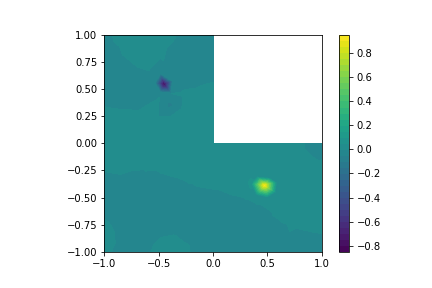}
    \caption{\rem{Weighted sparsity regularization with $\alpha=10^{-3}$ recovers the true sink-source configuration shown in figure \ref{fig:motivation}(a).}}
    \label{fig:motivation2}
\end{figure}

\rem{
Recall that standard sparsity regularization did not produce satisfactory results for the example considered in section \ref{sec:motivation}, see figure \ref{fig:motivation}. Figure \ref{fig:motivation2} shows that the weighted version \eqref{eq:5} handles this case very well.
}

\subsection{Example 1: All but one in $\NAT$}
\begin{figure}[H]
    \centering
    \begin{subfigure}[b]{0.5\linewidth}        
        \centering
        \includegraphics[width=\linewidth]{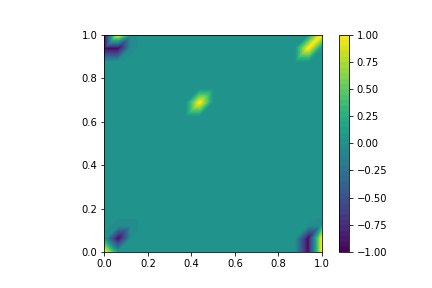}
        \caption{True sinks and sources.}
    \end{subfigure}\par
    \begin{subfigure}[b]{0.5\linewidth}        
        \centering
        \includegraphics[width=\linewidth]{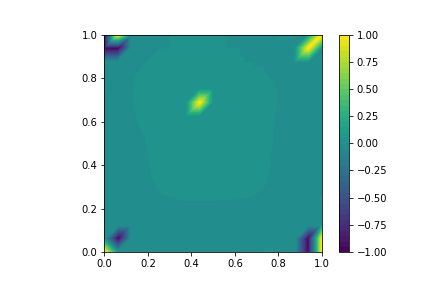}
        \caption{The solution $\yy_\alpha^*$ of \eqref{eq:6}, i.e., using weighted regularization.}
    \end{subfigure}\par
    \begin{subfigure}[b]{0.5\linewidth}        
        \centering
        \includegraphics[width=\linewidth]{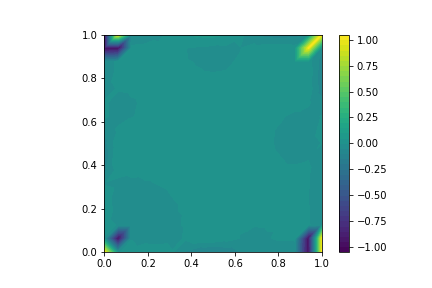}
        \caption{Inverse solution computed with unweighted/standard regularization.}
    \end{subfigure}\par
    \caption{Example 1. Comparison of the true sinks and sources and the inverse solutions. The regularization parameter was $\alpha = 10^{-4}$.}
    \label{fig:ex1}
\end{figure}

The \rem{second} example concerns Corollary \ref{reg_Th_boundary_plusOne}: All but one of the sources or sinks are in the orthogonal complement $\NAT$ of the null space $\nullspace{\AAA}$ of the forward matrix $\AAA$. To set up the experiment, we selected a number of basis vectors $\{ \ee_j \}_{j \in \mathcal{J}\setminus\{\notj\}}$ satisfying $\|\PP\ee_j\|_2 \geq 0.95$, $j \in \mathcal{J}\setminus\{\notj\}$, and arbitrarily set them to be sinks or sources with unit magnitude. Here, we recall that $\|\PP\ee_j\|_2 = 1$ implies that $\ee_j \in \NAT$, see \eqref{eq:3}. Thereafter we added one interior source, i.e., $\ee_{\notj}$. The outcome of this process is illustrated in figure \ref{fig:ex1}(a), which shows the true sources and sinks. 

Panels (b) and (c) in figure \ref{fig:ex1} display the regularized solutions when weighted and unweighted regularization are employed, respectively. We observe that all the sinks and sources located at the boundary are recovered with both methods, but the interior source is only detected when weighted regularization is employed.   

\subsection{Example 2: Well-separated sources and sinks} 
According to Corollary \ref{reg_Th_disjoint_projections}, a collection of well-separated sinks and sources can be recovered. That is, this corollary requires that the projections $\{ \PP \ee_j \}_{j \in \mathcal{J}}$ 
are disjoint. 

\begin{figure}[H]
    \centering
    \begin{subfigure}[b]{0.49\linewidth}        
        \centering
        \includegraphics[width=\linewidth]{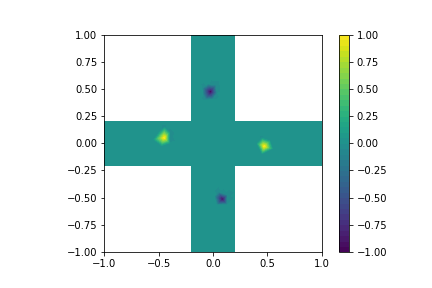}
        \caption{True sinks and sources.}
    \end{subfigure}
    \begin{subfigure}[b]{0.49\linewidth}        
        \centering
        \includegraphics[width=\linewidth]{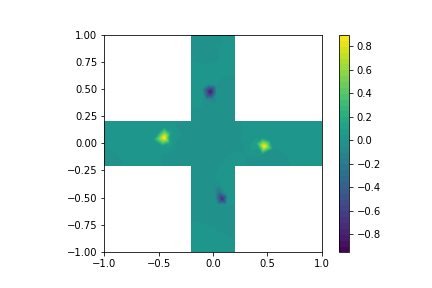}
        \caption{0\% noise ($\alpha = 10^{-4}$).}
    \end{subfigure}\par
    \begin{subfigure}[b]{0.49\linewidth}        
        \centering
        \includegraphics[width=\linewidth]{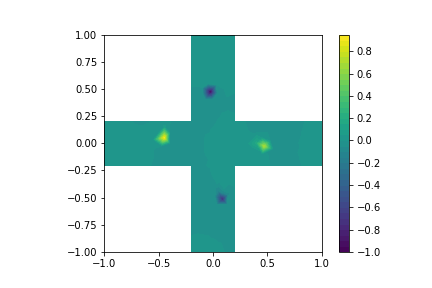}
        \caption{1\% noise ($\alpha = 0.005$).}
    \end{subfigure}
    \begin{subfigure}[b]{0.49\linewidth}        
        \centering
        \includegraphics[width=\linewidth]{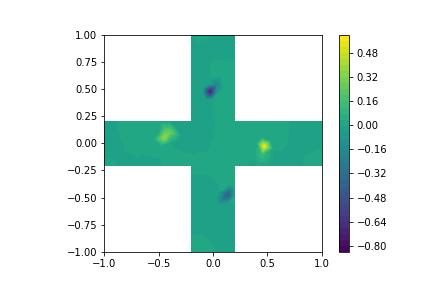}
        \caption{5\% noise ($\alpha = 0.025$).}
    \end{subfigure}    
    \caption{Example 2. Comparison of the true sinks and sources and the inverse solutions computed with weighted regularization and different levels of noise.}
    \label{fig:ex2.1}
\end{figure}

\begin{figure}[H]
    \centering
    \begin{subfigure}[b]{0.49\linewidth}        
        \centering
        \includegraphics[width=\linewidth]{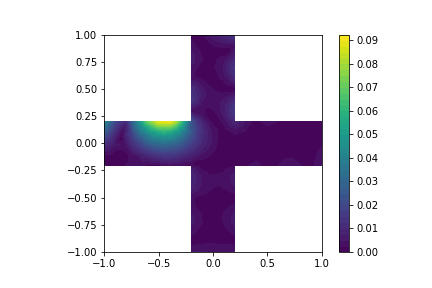}
    \end{subfigure}
    \begin{subfigure}[b]{0.49\linewidth}        
        \centering
        \includegraphics[width=\linewidth]{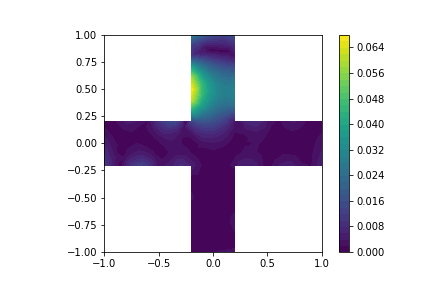}
    \end{subfigure}\par
    \begin{subfigure}[b]{0.49\linewidth}        
        \centering
        \includegraphics[width=\linewidth]{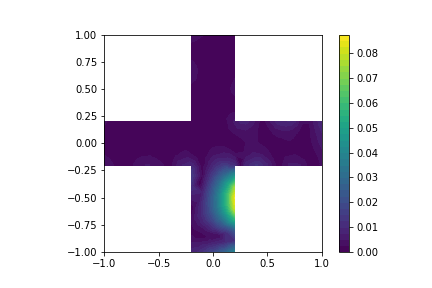}
    \end{subfigure}
    \begin{subfigure}[b]{0.49\linewidth}        
        \centering
        \includegraphics[width=\linewidth]{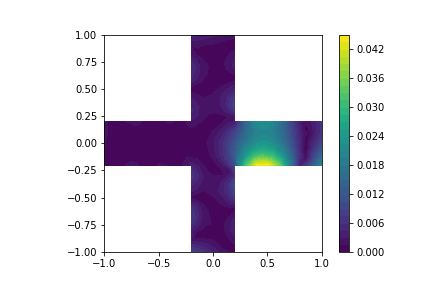}
    \end{subfigure}    
    \caption{Example 2. Plots of $|\PP\ee_j|$ for the four individual sources/sinks displayed in Figure 2a).}
    \label{fig:ex2.4}
\end{figure}

Panel (a) in figure \ref{fig:ex2.1} visualises the cross-shaped domain and the true sink-source configuration considered in our \rem{third} test problem. Figure \ref{fig:ex2.4} shows that the supports of $\{ \PP \ee_j \}_{j \in \mathcal{J}}$ are almost disjoint. 
(In figure \ref{fig:ex2.4} we visualise, for the sake of easy interpretation, the vector $|\PP \ee_j|$ containing the absolute values of the components of $\PP \ee_j$.)

The source-sink detection works well in this case, see panels (b), (c) and (d) in figure \ref{fig:ex2.1}. More precisely, the positions of the individual sinks and sources are perfectly recovered, but the magnitude is underestimated when the noise level increases. We also mention that no source is misinterpreted as a sink, or vice versa. 

Concerning the noise, we computed the synthetic observation data 
\begin{equation*}
   \bb = \AAA\xx^* + \tau\boldsymbol{\rho},
\end{equation*}
where $\tau$ is a scalar, $\boldsymbol{\rho}$ is a vector containing normally distributed numbers with zero mean and standard deviation $1$, and $\xx^*$ represents the true sources and sinks. The noise level is then defined as the ratio 
\begin{equation*}
    \frac{\tau}{\max\{\bb\} - \min\{\bb\}}. 
\end{equation*}
We used Morozov's discrepancy principle to select the size of the regularization parameter $\alpha$. 

\begin{figure}[H]
    \centering
        \includegraphics[width=0.7\linewidth]{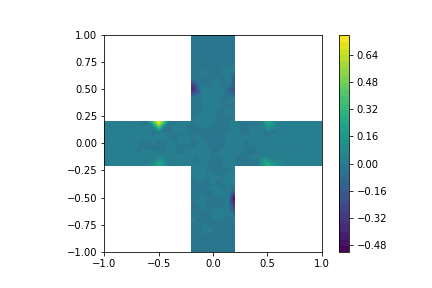}
    \caption{Example 2. The inverse solution computed with unweighted regularization. The regularization parameter was $\alpha = 10^{-4}$. Panel (a) in figure \ref{fig:ex2.1} shows the true sinks and sources.}
    \label{fig:ex2.2}
\end{figure}

\begin{figure}[H]
    \centering
    \begin{subfigure}[b]{0.49\linewidth}        
        \centering
        \includegraphics[width=\linewidth]{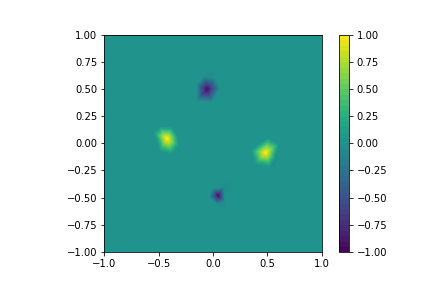}
        \caption{True sinks and sources.}
    \end{subfigure}
    \begin{subfigure}[b]{0.49\linewidth}        
        \centering
        \includegraphics[width=\linewidth]{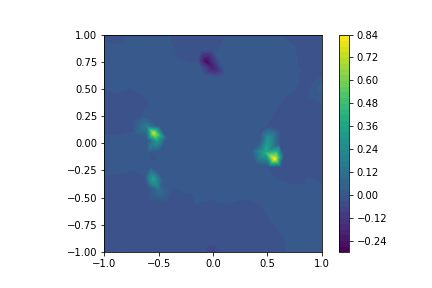}
        \caption{Inverse solution.}
    \end{subfigure}
    \caption{Example 2. Comparison of the true sources and sinks and the inverse solution computed with weighted regularization. The regularization parameter was $\alpha = 10^{-4}$.}
    \label{fig:ex2.3}
\end{figure}

According to figure \ref{fig:ex2.2}, standard unweighted regularization fails to handle this case adequately. We also observe that the results generated with the weighting procedure deteriorates when a square domain $\Omega=(-1,1)^2$ is used instead of a cross-shaped domain, see figure \ref{fig:ex2.3}. In the square domain case, the supports of $\{ \PP \ee_j \}_{j \in \mathcal{J}}$ become much more overlapping (illustration not included) than with the cross-shaped geometry, cf. Corollary \ref{reg_Th_disjoint_projections}. Roughly speaking, source-sink detection is more difficult for convex domains than for non-convex geometries.  

\begin{figure}[H]
    \centering
    \begin{subfigure}[b]{0.49\linewidth}        
        \centering
        \includegraphics[width=\linewidth]{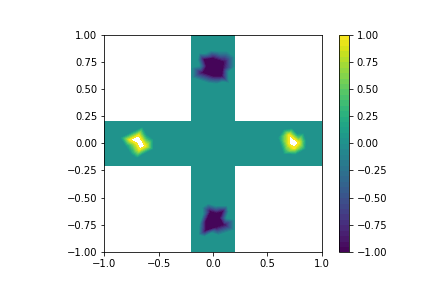}
        \caption{True sinks and sources.}
    \end{subfigure}
    \begin{subfigure}[b]{0.49\linewidth}        
        \centering
        \includegraphics[width=\linewidth]{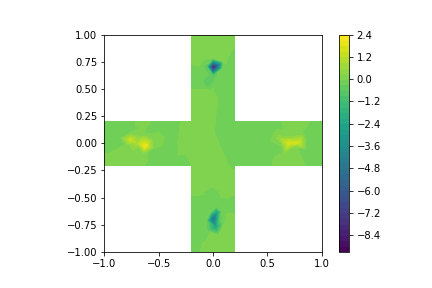}
        \caption{Inverse solution.}
    \end{subfigure}
    \caption{Example 2. Comparison of the true composite sources and sinks and the inverse solution computed with weighted regularization. The regularization parameter was $\alpha = 10^{-4}$.}
    \label{fig:ex2.5}
\end{figure}

Figure \ref{fig:ex2.5} shows that the weighting procedure might not work very well for determining the size of composite sinks and sources. That is, sources and sinks that are not generated by a single basis vector $\ee_j$. The inverse solution correctly identifies the positions of the sources and sinks, but the magnitudes and extends are not correct. How to handle this, is an open problem. Will box constraints resolve this issue?  

\rem{
\subsection{Example 3: Convergence rates}
We end this section with a numerical investigation inspired by Theorem \ref{thm:convergence}. The true solution $\xx^*$ in this case consists of one source and one sink, and we use the same cross-shaped domain as in the previous example. We present two loglog-plots of how the errors \eqref{eq:convX} and \eqref{eq:convY} decrease as the noise level $\delta$ decreases. That is, we consider the regularized problems \eqref{eq:formA} and \eqref{eq:formAd}.
}

\begin{figure}[H]
    \centering
    \begin{subfigure}[b]{1.1\linewidth}        
        \centering
        \includegraphics[width=0.6\linewidth]{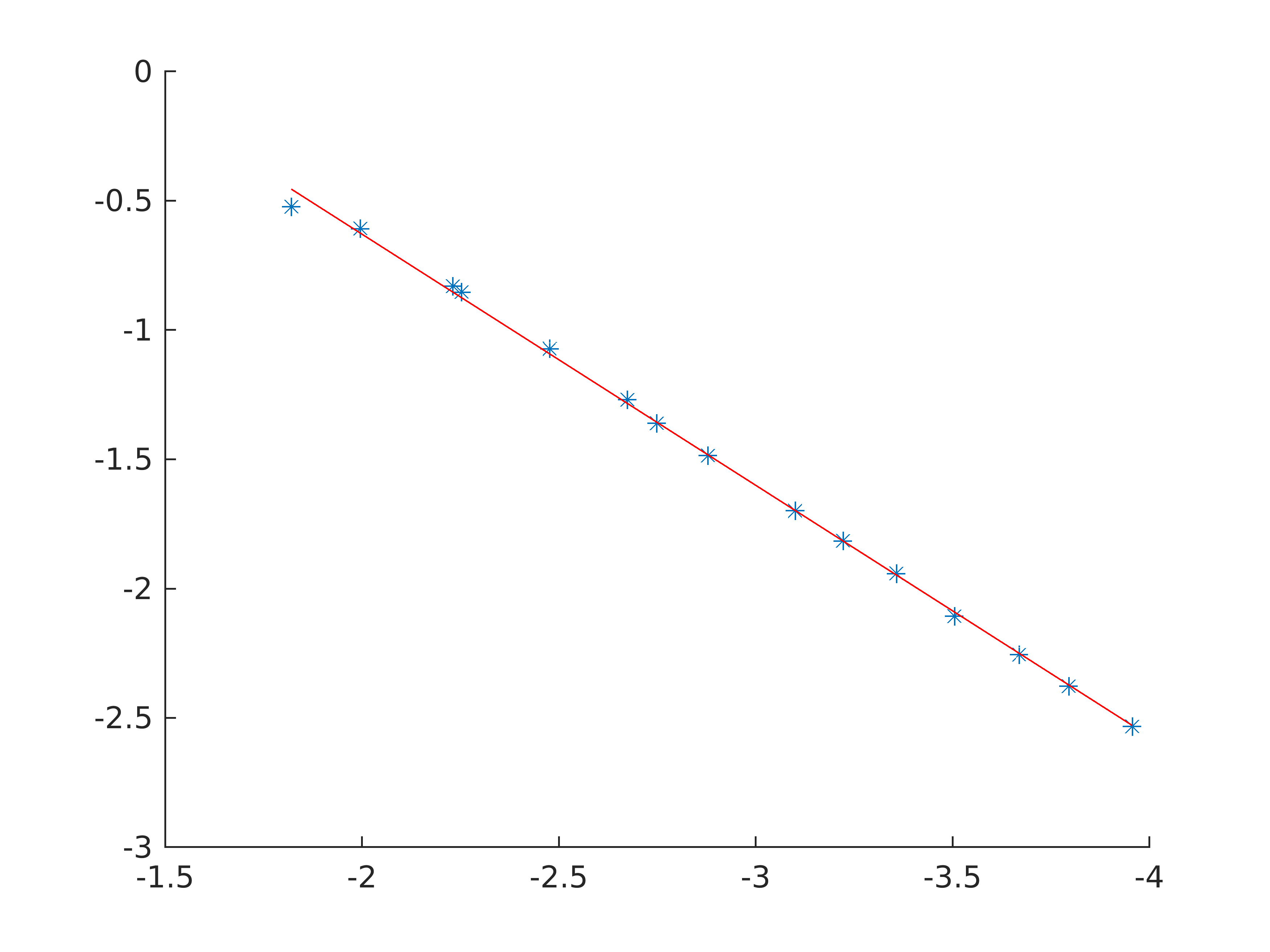}
        \caption{Error term $\|\xx_\alpha - \xx^*\|_{\WW}$, cf. \eqref{eq:convX}}
    \end{subfigure}\par
    \begin{subfigure}[b]{1.1\linewidth}        
        \centering
        \includegraphics[width=0.6\linewidth]{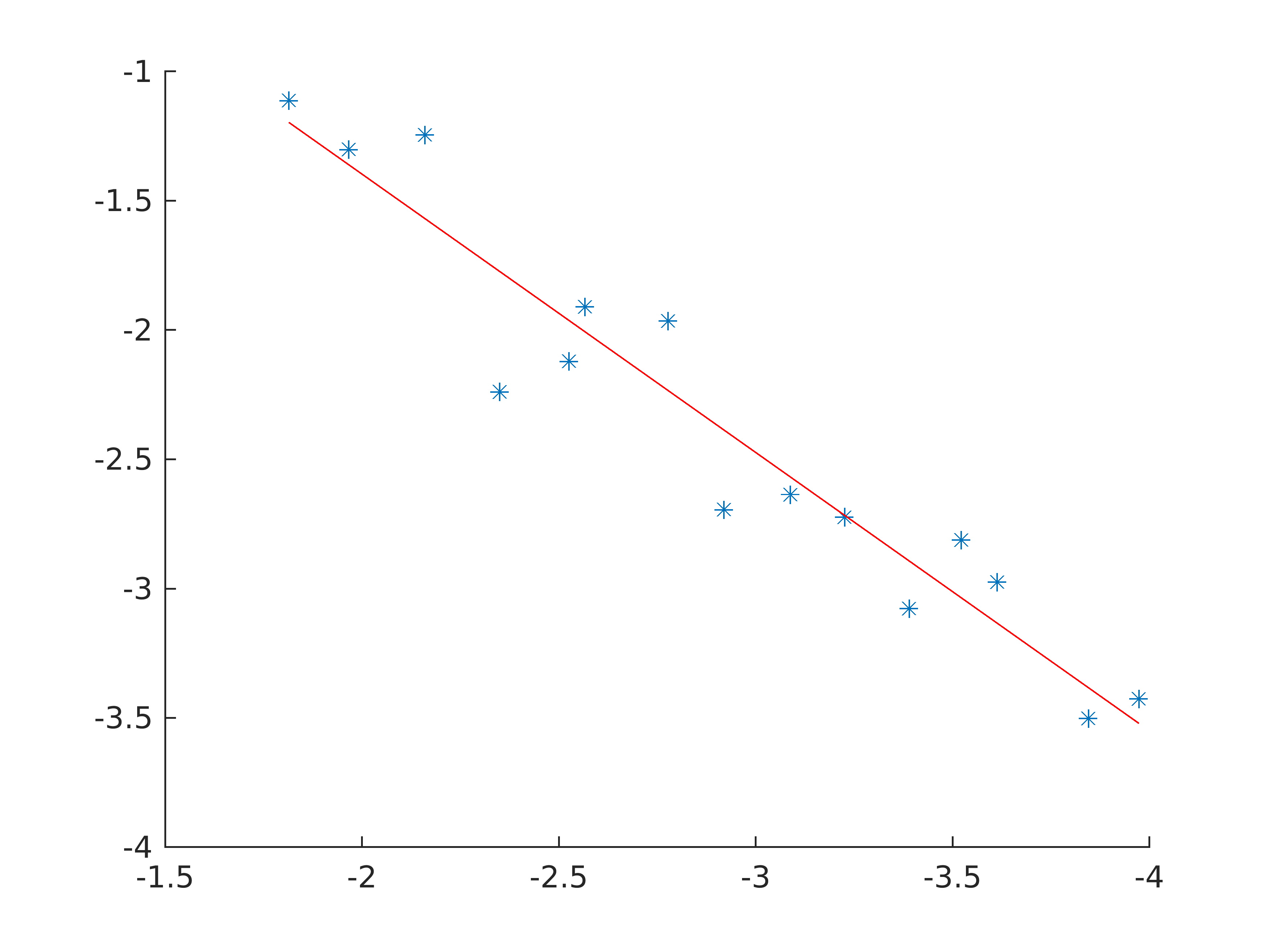}
        \caption{Error term $\|\yy_\alpha - \xx^*\|_{\WW}$, cf. \eqref{eq:convY}}
    \end{subfigure}
    \caption{Example 3. \rem{Loglog-plots of the noise level $\delta$ (horizontal axis) versus the error term $\|\qq_\alpha - \xx^*\|_{\WW}$ (vertical axis) for $\qq_\alpha$ equal to $\xx_\alpha$ and $\yy_\alpha$. The blue asterisks represent the numerical results, i.e., the error $\|\qq_\alpha - \xx^*\|_{\WW}$,  and the red lines are the outcome of employing linear regression to the log-values of the data points. }}
    \label{fig:ex3}
\end{figure}

\rem{
In figure \ref{fig:ex3}(a) we clearly observe the linear convergence behaviour \eqref{eq:convX} with only small deviations from a linear trend. When we instead consider the formulation \eqref{eq:formAd}, we still observe a linear trend, but with much larger deviations from a "regular" linear behaviour, see figure \ref{fig:ex3}(b). It is reasonable to believe that the latter phenomenon is due to the fact that the noise vector is multiplied by $\AAA_k^\dagger$. 
}

\appendix 
\section{\rem{Proof of Theorem \ref{thm_main_basis_pursuit}}} \label{app:proof}
\begin{proof}
Let $$\yy = \sum_{\JJ^c \cup \JJ} y_i \ee_i$$ denote a solution of \eqref{NBP3}. Then 
\[
\AAA \yy = \AAA \sum_{\JJ} x_j^* \ee_j
\]
or, because $\PP = \AAA^\dagger \AAA$,  
\[
\PP \yy = \PP \sum_{\JJ} x_j^* \ee_j. 
\]
Taking the inner product with $\cc$ yields 
\begin{equation} \label{NBP3.1}
\sum_{\JJ^c \cup \JJ} y_i \PP \ee_i \cdot \cc  
= \sum_{\JJ} x_j^* \PP \ee_j \cdot \cc 
\end{equation}
where, using \eqref{NBP1} and \eqref{eq:4},  
\begin{align}  
\nonumber
\sum_{\JJ} x_j^* \PP \ee_j \cdot \cc &= \sum_{\JJ} x_j^* \, \textnormal{sgn}(x_j^*) \,  \| \PP \ee_j \|_2 \\
\nonumber
&= \sum_{\JJ} |x_j^*| \, \| \PP \ee_j \|_2 \\ 
\label{NBP4}
&= \| \WW \xx^* \|_1
\end{align}
and, invoking \eqref{NBP1}, \eqref{NBP2} and \eqref{eq:4}, 
\begin{align}
\nonumber
\sum_{\JJ^c \cup \JJ} y_i \PP \ee_i \cdot \cc &\leq \sum_{\JJ^c \cup \JJ} | y_i| \left| \PP \ee_i \cdot \cc \right| \\
\nonumber
&\leq \sum_{\JJ^c \cup \JJ} | y_i| \left\| \PP \ee_i \right\|_2 \\
\label{NBP4.01}
&= \| \WW \yy \|_1.
\end{align}
From \eqref{NBP3.1}, \eqref{NBP4} and \eqref{NBP4.01} we can conclude that 
\begin{equation*}
\| \WW \yy \|_1 \geq \| \WW \xx^* \|_1    
\end{equation*}
and it follows that $\xx^*$ is a solution of \eqref{NBP3}.

If there is $i' \in \JJ^c$ such that $y_{i'} \neq 0$, then the strict inequality in \eqref{NBP2} yields a strict inequality in \eqref{NBP4.01}:  
\begin{align*}
\sum_{\JJ^c \cup \JJ} y_i \PP \ee_i \cdot \cc &\leq \sum_{\JJ} | y_i| \left| \PP \ee_i \cdot \cc \right| + \sum_{\JJ^c} | y_i| \left| \PP \ee_i \cdot \cc \right| \\
&< \sum_{\JJ} | y_i| \left\| \PP \ee_i \right\|_2 + \sum_{\JJ^c} | y_i| \left\| \PP \ee_i \right\|_2 \\
&= \| \WW \yy \|_1,
\end{align*}
and this would imply that $\| \WW \yy \|_1 > \| \WW \xx^* \|_1$. We conclude that any solution $\yy$ of \eqref{NBP3} must satisfy $\textnormal{supp}(\yy) \subseteq \textnormal{supp}(\xx^*)$.
\end{proof}

\bibliographystyle{abbrv}
\bibliography{references}

\end{document}